\documentstyle[12pt]{article}     
\begin{document} 
\title{Stochastic processes on non-Archimedean spaces. II.
Stochastic antiderivational equations.}  
\author{S.V. Ludkovsky}  
\date{23 March 2001
\thanks{Mathematics subject classification
(1991 Revision) 28C20 and 46S10.} }
\maketitle
\par address: Laboratoire de Math\`ematiques Pures, 
\par Complexe Scientifique des C\`ezeaux,
\par 63177 AUBI\`ERE Cedex, France.
\par permanent address: Theoretical Department,
Institute of General Physics,
Str. Vavilov 38, Moscow, 119991 GSP-1, Russia.
\begin{abstract} 
Stochastic antiderivational equations on Banach spaces over 
local non-Archimedean fields are investigated.
Theorems about existence and uniqueness of the
solutions are proved under definite conditions.
In particular Wiener processes are considered
in relation with the non-Archimedean analog of the Gaussian measure.
\end{abstract} 
\section{Introduction.}
This article continues investigations of stochastic processes on
non-Archimedean spaces (\cite{lust1}).
In the first part stochastic processes were defined
on Banach spaces over non-Archimedean local fields 
and the analogs of It$\hat o$ formula
were proved. This part is devoted to stochastic antiderivational equations.
In the non-Archimedean case antiderivational equations are used instead of
stochastic integral or differential equations in the classical case.
\par Stochastic differential equations on real Banach spaces and manifolds
are widely used for solutions of mathematical and physical
problems and for construction and investigation of measures
on them \cite{dalf,gihsko,ikwat,mckean,malli,oksen}.
Wide classes of quasi-invariant measures including 
analogous to Gaussian type
on non-Archimedean Banach spaces, loops and diffeomorphisms 
groups were investigated in \cite{lud,lu6,luumnls,lutmf99,lubp2}.
Quasi-invariant measures on topological groups and their
configuration spaces can be used 
for the investigations of their unitary representations
(see \cite{luseamb,lubp99,lutmf99,lubp2} and references therein).
In view of this developments non-Archimedean analogs
of stochastic equations and diffusion processes 
need to be investigated. Some steps in this direction were 
made in \cite{bikvol,evans}. 
There are different variants for such activity, for example,
$p$-adic parameters analogous to time, but spaces
of complex-valued functions. At the same time measures may be 
real, complex or with values in a non-Archimedean field.
In the classical stochastic analysis indefinite integrals
are widely used, but in the non-Archimedean case they have quite another
meaning, because the field of $p$-adic
numbers $\bf Q_p$ has not any linear order structure compatible
with its normed field structure (see Part I).
\par This work treats the case which was not considered by 
another authors and that is suitable and helpful for 
the investigation of stochastic processes and quasi-invariant measures
on non-Archimedean topological groups. 
In \S 2 suitable analogs of Gaussian measures are considered.
Certainly they have not any complete analogy with the classical one,
some of their properties are similar and some are different.
They are used for the definiton of the standard (Wiener)
stochastic process.
Integration by parts formula for the non-Archimedean stochastic processes
is studied. Some particular cases of the general It$\hat o$ formula
from Part I are dicussed here more concretely.
In \S 3 with the help of them stochastic antiderivational equations are 
defined and investigated. Analogs of 
theorems about existence and uniquiness of solutions
of stochastic antiderivational equations are proved.
Generating operators of solutions of stochastic equations are investigated.
All results of this paper are obtained for the first time.
\par In this part the notations of Part I also are used.  
\section{Gaussian measures and standard Wiener processes
on a non-Archimedean Banach space.}
\par {\bf 2.1.} Let $H=c_0(\alpha ,{\bf K})$ be a Banach space over
a local field $\bf K$ with an ordinal $\alpha $
and the standard orthonormal base $\{ e_j: j\in \alpha \} $,
$e_j=(0,...,0,1,0,...)$ with $1$ on the $j$-th place.
Let $\sf U^P$ be a cylindrical algebra generated by projections on
finite-dimensional over $\bf K$ subspaces $F$ in $H$ and Borel 
$\sigma $-algebras $Bf(F)$. Denote by $\sf U$ the minimal $\sigma $-algebra
$\sigma ({\sf U^P})$ generated by $\sf U^P$. When
$card (\alpha )\le \aleph _0$, then ${\sf U}=Bf(H)$.
Each vector $x\in c_0$ is considered as continuous linear functional 
on $c_0$ by the formula $x(y)=\sum_jx^jy^j$ for each $y\in c_0$, 
so there is the natural embedding $c_0\hookrightarrow c_0^*$,
where $x=\sum_jx^je_j$, $x^j\in \bf K$.
The field $\bf K$ is the finite algebraic extension of
of the field $\bf Q_p$ of $p$-adic numbers and as the Banach space over
$\bf Q_p$ it is isomorphic with $\bf Q_p^n$, that is, each
$z\in \bf K$ has the form $z=(z^1,...,z^n)$, where $z^1,...,z^n\in \bf Q_p$.
Let $\{ y \} _p :=\sum_{j<0}y_jp^j$, where $y\in \bf Q_p$, $y=\sum_jy_jp^j$,
$y_j \in \{ 0,1,...,p-1 \} $, in particular for values
$y=(z,x):=\sum_{j=1}^nx^jz^j$ for $x, z\in \bf K$.
All continuous characters $\chi : {\bf K}\to \bf C$ of $\bf K$ as
the additive group have the form 
$$(i)\mbox{ }\chi _{\gamma }(x)=\epsilon ^{z^{-1}\{ (e,\gamma x) \} _p}$$ 
for each $ \{ (e,\gamma x) \} _p \ne 0$,
$\chi _{\gamma }(x):=1$ for $ \{ (e,\gamma x) \} _p=0,$
where $\epsilon =1^z$ is a root of unity, $z=p^{ord(\{ (e,\gamma x) \} _p)},$ 
$e=(1,...,1)\in \bf Q_p^n$, $\gamma \in \bf K$ (see \S 25 \cite{hew}
and \S I.3.6, about the spaces $L_q(H)$ of operators see \S I.2). 
Each $\chi $ is locally constant, hence $\chi : {\bf K}\to
\bf T$ is also continuous, where $\bf T$ denotes
the discrete group of all roots of $1$ (by multiplication).
\par Let us consider functions, whose Fourier transform has the form:
$$\hat f(x)=\hat f_{\beta ,\gamma ,q}(x):=exp
(-\beta |x|^q)\chi _{\gamma }(x),$$
where the Fourier transform was defined in \S 7
\cite{vla3} and \cite{roo}, $\gamma \in \bf K$, $0<\beta <\infty $,
$0<q<\infty $.
\par {\bf Definition.} A cylindrical measure $\mu $ on $\sf U^P$
is called $q$-Gaussian, if each its one-dimensional projection 
is $q$-Gaussian, that is, 
$$(i)\quad \mu ^g(dx)=C_{\beta ,\gamma ,q}
f_{\beta ,\gamma ,q}v(dx),$$ 
where $v$ is the Haar measure on $Bf({\bf K})$ with values 
in $\bf R$, where $g$ is a continuous $\bf K$-linear 
functional on $H=c_0(\alpha ,{\bf K})$ 
giving projection on one-dimensional subspace in $H$,
$C_{\beta ,\gamma ,q}>0$ are constants such that $\mu ^g({\bf K})=1$, 
$\beta $ and $\gamma $ may depend on $g$, $q$ is independent of $g$,
$1\le q<\infty $, $\alpha \subset \omega _0$, $\omega _0$ is the first
countable ordinal.
\par If $\mu $ is a measure on $H$, then
$\hat \mu $ denotes its characteristic functional, that is,
$\hat \mu (g):=\int_H\chi _g(x)\mu (dx)$, where
$g\in H^*$, $\chi _g: H\to \bf C$ is the character
of $H$ as the additive group (see \S I.3.6).
\par {\bf 2.2.} {\bf Theorem.} {\it A non-negative $q$-Gaussian measure 
$\mu $ on $c_0(\omega _0,{\bf K})$ is $\sigma $-additive on 
$Bf(c_0(\omega _0,{\bf K}))$ if and only if
there exists an injective compact operator 
$J\in L_q(c_0(\omega _0,{\bf K}))$ for a chosen $1\le q<\infty $
such that 
$$(i)\mbox{ }\mu (dx)=
\bigotimes_{j=1}^{\infty }\mu _j(dx^j),\mbox{ }where $$
$$(ii)\mbox{ }J=diag(\zeta _j: \mbox{ }\zeta _j\in {\bf K},
\mbox{ }j\in \omega _0),$$
$$(iii)\mbox{ }\mu _j(dx^j)=C_{\beta _j,\gamma _j,q}f_{\beta _j,\gamma _j,q}
v(dx^j)$$
are measures on $e_j\bf K$, $x=(x^j: j\in \omega _0)
\in c_0(\omega _0,{\bf K})$, $x^j\in \bf K$, $\beta _j=|\zeta _j|^{-q}$,
$\gamma =(\gamma _j: j\in \omega _0)\in c_0(\omega _0,{\bf K})$.
Moreover, each one-dimensional projection $\mu ^g$
has the following characteristic functional:
$$(iv)\quad {\hat \mu }^g(h)=exp(-(\sum_j\beta _j|g_j|^q)|h|^q)
\chi _{g(\gamma )}(h),$$
 where $g=(g_j: j\in \omega _0)\in c_0(\omega _0,{\bf K})^*$.}
\par {\bf Proof.} Let $\theta $ be a characteristic functional
of $\mu $. By the non-Archimedean analog of the Minlos-Sazonov
Theorem  (see \S 2.31 in \cite{lu6}, \cite{lubkt}) a measure $\mu $ is 
$\sigma $-additive if and only if for each $c>0$ there exists
a compact operator $S_c$ such that $|Re(\mu (y)-\mu (x))|
<c$ for each $x, y\in c_0(\omega _0,{\bf K})$ with $|z^*S_cz|<1$, where
$z=x-y$.
From the definition of $\mu $ to be $q$-Gaussian it follows, that
each its projection $\mu _j$ on ${\bf K}e_j$ has the form given by Equation
$(iii)$. It remains to establish that $\mu $ is $\sigma $-additive 
if and only if $J\in L_q(c_0(\omega _0,{\bf K}))$ and 
$\gamma \in c_0(\omega _0,{\bf K})$.
\par In view of Lemma 2.3 \cite{lu6} $\mu $ is $\sigma $-additive
if and only if each sequence of finite-dimensional (over $\bf K$ 
distributions) satisfies two conditions:\\
$(2.3.i)$ for each $c>0$ there exists $b>0$ such that
$\sup_n|\mbox{ }|\mu _{L(n)}|(B(c_0,0,r)\cap L(n))-
|\mu _{L(n)}|(L(n))|\le c$ for each $r\ge b$,\\
$(2.3.ii)$ $\sup_n|\mu _{L(n)}|(L(n))<\infty $.
Take in particular $L(n)=sp_{\bf K}\{ e_1,...,e_n \} $
for each $n\in \bf N$.
\par We have $\mu _j({\bf K}\setminus B({\bf K},0,r))
\le C\mbox{ }\int_{x\in {\bf K},|x|>r}exp(-|x/\zeta _j|^q)|\zeta _j|^{-1}v(dx)$
$\le C_1\int_{y\in {\bf R},|y|>r}exp(-|y|^q|\zeta _j|^{-q})|\zeta _j|^{-1}dy,$ 
where $C>0$ and $C_1>0$ are constants independent from 
$\zeta _j$ for $b_0>p^3$ and each $r>b_0$, $1\le q<\infty $ is fixed
(see also the proof of Lemma 2.8 \cite{lu6}
and Theorem II.2.1 \cite{dalf}). Evidently, $g(\gamma )$ 
is correctly defined for each $g\in c_0(\omega _0,{\bf K})^*$
if and only if $\gamma \in c_0(\omega _0,{\bf K})$. In this case
the character $\chi _{g(\gamma )}: {\bf K}\to \bf C$ is defined
and $\chi _{g(\gamma )}=\prod_{j=1}^{\infty }\chi _{g_j\gamma _j}$.
Therefore, if $J\in L_q(c_0)$ and $\gamma \in c_0(\omega _0,{\bf K})$, then
$\mu $ is $\sigma $-additive.
\par Let $0\ne g\in c_0^*$. Since $\bf K$ is the local field
there exists $x_0\in c_0$ such that $|g(x_0)|=\| g\| $
and $\| x_0 \| =1.$ Put $g_j:=g(e_j)$. Then $\| g\|
\le \sup_j|g_j|$, since $g(x)=\sum_jx^jg_j$, where $x=x^je_j:=\sum_jx^je_j$
with $x^j\in \bf K$. Consequently, $\| g\| =\sup_j|g_j|$.
We denumerate the standard orthonormal basis
$\{ e_j: j\in {\bf N} \} $ such that $|g_1|=\| g\| $.
There exists an operator $E$ on $c_0$ with matrix elements
$E_{i,j}=\delta _{i,j}$ for each $i, j>1$,
$E_{1,j}=g_j$ for each $j\in {\bf N}$.
Then $|det\mbox{ }P_nEP_n|=\| g \|$ for each $n\in \bf N$, where
$P_n$ are the standard projectors on $sp_{\bf K}\{ e_1,...,e_n \} $.
When $g\in \{ e_j^*: j\in \omega _0 \} $, then evidently,
$\mu ^g$ has the form given by Equation $(iii)$,
since $\mu _i({\bf K})=1$ for each $i\in \omega _0 $, where
$e_j^*(e_i)=\delta _{i,j}$ for each $i, j$.
\par Suppose now that $J\notin L_q(c_0)$. 
For this we consider $\mu ^g({\bf K}\setminus B({\bf K},0,r))\ge 
\sum_j\int_{x\in {\bf K},|x|>r}C\mbox{ }
exp(-|x/\zeta _j|^q)|\zeta _j|^{-1}v(dx)$,
where $g=(1,1,1,...)
\in c_0^*=l^{\infty }(\omega _0,{\bf K})$.
On the other hand, there exists a constant $C_2>0$ such that for 
$b_0>p^3$ and each $r>b_0$ there is the following inequality:
$\int_{x\in {\bf K}, |x|>r}C\mbox{ }
exp(-|x/\zeta _j|^q)|\zeta _j|^{-1}v(dx) \ge $ \\
$C_2[\int_r^{\infty }exp(-|y|^q|\zeta _j|^{-q})|\zeta _j|^{-1}dy
+ \int_{-\infty }^{-r}exp(-|y|^q|\zeta _j|^{-q})|\zeta _j|^{-1}dy]$. 
From the estimates 
of Lemma II.1.1 \cite{dalf} and using the substitution
$z=y^{1/2q}$ for $y>0$ and $z=(-y)^{1/2q}$ for $y<0$
we get that $\mu ^g$
is not $\sigma $-additive, consequently, $\mu $ is not 
$\sigma $-additive, since $P_g^{-1}(A)$ are cylindrical Borel 
subsets for each $A\in Bf({\bf K})$, where $P_gz=g(z)$ is 
the induced projection on $\bf K$ for each $z\in c_0$. 
\par For the verification of Formula $(iv)$ it is sufficient at first to
consider the measure $\mu $ on the algebra $\sf U^P$ of cylindrical subsets
in $c_0$. Then for each projection $\mu ^g$, where $g\in sp_{\bf K}
(e_1,...,e_m)^*$, we have:
$${\hat \mu }^g(h)=\int_{\bf K}...\int_{\bf K}\chi _e(hz)
\mu _1(dx_1)...\mu _m(dx_m),$$ 
where $e=(1,...,1)\in \bf Q_p^n$, $h\in
{\bf K}$, $n:=dim_{\bf Q_p}{\bf K},$
$x^i\in {\bf K}e_i$, $z=g(x)$, $x=(x^1,...,x^m)$, consequently,
$\hat \mu ^g(h)=\prod_{i=1}^m\hat \mu _i(hg_i)$, since
$\chi _e(hg(x))=\prod_{i=1}^m\chi _e(h_ig_ix^i)$
for each $x\in sp_{\bf K}(e_1,...,e_m).$
Since $J\in L_q$, then $\mu $ is the Radon measure, consequently,
the continuation of $\mu $ from $\sf U^P$ produces
$\mu $ on the Borel $\sigma $-algebra of $c_0$,
hence $\lim_{m\to \infty }{\hat \mu }^{Q_mg}(h)=
\hat \mu ^g(h),$ where $Q_m$ is the natural projection on
$sp_{\bf K}(e_1,...,e_m)^*$ for each $m\in \bf N$
such that $Q_m(g)=(g_1,...,g_m)$. Using expressions
of $\hat \mu _i$ we get Formula $(iv)$.
From this follows, that if $J\in L_q$, then $\hat \mu (g)$ exists for each
$g\in c_0^*$ if and only if $\gamma \in c_0$, since $\hat \mu ^g(h)=
\hat \mu (gh)$ for each $h\in \bf K$ and $g\in c_0^*$.
\par {\bf 2.3. Corollary.} {\it $|\hat \mu ^g(h_1+h_2)|
\le \max (|\hat \mu ^g(h_1)|, |\hat \mu ^g(h_2)|)$
for each $h_1, h_2 \in \bf K$ and $g\in c_0(\omega _0,{\bf K})^*$.}
\par {\bf Proof.} In view of the ultrametric inequality
$|h_1+h_2|^q\le \max (|h_1|^q,|h_2|^q)$ for each $1\le q<\infty $
and $h_1, h_2 \in \bf K$. Since $|\chi _{\gamma }(h)|=1$ for each
$h, \gamma \in \bf K$, then from
Formula 2.2.(iv) the statement of this Corollary follows.
\par {\bf 2.4. Remark.} Let $Z$ be a compact subset without
isolated points in a local field $\bf K$, for example,
$Z=B({\bf K},t_0,1)$. Then 
the Banach space $C^0(Z,{\bf K})$ has the Amice polynomial 
orthonormal base
$Q_m(x)$, where $x\in Z$, $m\in {\bf N_o}:=\{ 0,1,2,... \} $ \cite{ami}.
Suppose $\tilde P^{n-1}: C^{n-1}(Z,{\bf K})\to C^n(Z,{\bf K})$ 
are antiderivations
from \S 80 \cite{sch1}, where $n\in \bf N$. Each $f\in C^0$
has a decomposition $f(x)=\sum_ma_m(f)Q_m(x)$, where $a_m\in \bf K$.
These decompositions establish the isometric isomorphism
$\theta : C^0(Z,{\bf K})\to c_0(\omega _0,{\bf K})$
such that $\| f\|_{C^0}=\max_m|a_m(f)|=\| \theta (f)\|_{c_0}$.
Since $Z$ is homeomorphic with $\bf Z_p$, then 
$\tilde P^1\tilde P^0: C^0(Z,{\bf K})\to C^2(Z,{\bf K})$
is a linear injective compact operator such that $\tilde P^1\tilde P^0\in L_1$,
where $\tilde P^j$ here corresponds to 
$\tilde P_{j+1}: C^j\to C^{j+1}$ antiderivation 
operator by Schikhof (see also \S \S 54, 80 \cite{sch1} and \S I.2.1). 
The Banach space $C^2(Z,{\bf K})$ is dense
in $C^0(Z,{\bf K})$. Using Theorem 2.2 and Note I.2.3 for $q\ge 1$
we get a $q$-Gaussian 
measure on $C^0(Z,{\bf K})$, where $\tilde P^1\tilde P^0f=\sum_j\lambda _jP_jf$
and $Jf=\sum_j\zeta _jP_jf$ for each $f\in C^0$, we put
$|\lambda _j| |\pi |^q\le |\zeta _j|^q\le |\lambda _j|$ for each $j\in \bf N$,
$P_j$ are projectors, $\lambda _j, \zeta _j\in \bf K$, $p^{-1}\le |\pi |<1$,
$\pi \in \bf K$ and $|\pi |$ is the generator of the valuation group
of $\bf K$.
\par If $H=c_0(\omega _0,{\bf K}),$ then the Banach space
$C^0(Z,H)$ is isomorphic with the tensor product 
$C^0(Z,{\bf K})\otimes H$ (see \S 4.R \cite{roo}).
Therefore, the antiderivation $\tilde P^n$ on $C^n(Z,{\bf K})$
induces the antiderivation $\tilde P^n$ on $C^n(Z,H)$. 
If $J_i\in L_q(Y_i)$, then $J:=J_1\otimes J_2\in L_q(Y_1\otimes Y_2)$
(see also Theorem 4.33 \cite{roo}). Put 
$Y_1=C^0(Z,{\bf K})$ and $Y_2=H$, 
then each $J:=J_1\otimes J_2\in L_q(Y_1\otimes Y_2)$
induces the $q$-Gaussian measure $\mu $ on $C^0(Z,H)$
such that $\mu =\mu _1\otimes \mu _2$, where $\mu _i$ are $q$-Gaussian
measures on $Y_i$ induced by $J_i$ as above. In particular for $q=1$
we also can take $J_1=\tilde P^1\tilde P^0$. 
The $1$-Gaussian measure on $C^0(Z,H)$ 
induced by $J=J_1\otimes J_2\in L_1$ with 
$J_1=\tilde P^1\tilde P^0$ we call standard.
Analogously considering the following Banach subspace
$C^0_0(Z,H):=\{ f\in C^0(Z,H):$ $f(t_0)=0 \} $
and operators $J:=J_1\otimes J_2\in L_1(C^0_0(Z,{\bf K})\otimes H)$
we get the $1$-Gaussian measures $\mu $ on it also, where $t_0\in Z$
is a marked point. Certainly, we can take others operators
$J_1\in L_q(Y_1)$ not related with the antiderivation as above.
\section{Non-Archimedean stochastic antiderivational equations.}
\par {\bf 3.1.} A measurable space $(\Omega ,{\sf F})$ with a normalised 
non-negative measure $\lambda $ 
is called a probability space and is denoted by
$(\Omega ,{\sf F},\lambda )$, where $\sf F$ is a $\sigma $-algebra of 
$\Omega $. Points $\omega \in \Omega $ are called 
elementary events and values $\lambda (S)$  are called
probabilities of events $S\in \sf F$. A measurable map
$\xi : (\Omega ,{\sf F})\to (X,{\sf B})$ is called a random variable
with values in $X$, where $\sf B$ is a $\sigma $-algebra of $X$ (see \S I.4.1). 
\par {\bf 3.2.} We define a (non-Archimedean) Wiener process 
$w(t,\omega )$ with values in $H$ 
as a stochastic process such that:
\par $(i)$ the differences $w(t_4,\omega )-w(t_3,\omega )$ 
and $w(t_2,\omega )-w(t_1,\omega )$ are independent
for each chosen $\omega $, $(t_1,t_2)$ and $(t_3,t_4)$ with $t_1\ne t_2$,
$t_3\ne t_4$, either $t_1$ or $t_2$ is not in the two-element set
$ \{ t_3,t_4 \} ,$ where $\omega \in \Omega ;$
\par $(ii)$ the random variable $\omega (t,\omega )-\omega (u,\omega )$ has 
a distribution $\mu ^{F_{t,u}},$ where $\mu $ is a probability 
Gaussian measure on $C^0(T,H)$ described in \S \S 2.1, 2.4, 
$\mu ^g(A):=\mu (g^{-1}(A))$ for $g\in C^0(T,H)^*$ and each $A\in
Bf(C^0(T,H)),$ a continuous linear functional $F_{t,u}$ is given by
the formula $F_{t,u}(w):=w(t,\omega )-w(u,\omega )$ 
for each $w\in L^s(\Omega ,{\sf F},\lambda ;C^0_0(T,H)),$
where $1\le s\le \infty ;$ 
\par $(iii)$ we also put $w(0,\omega )=0,$  
that is, we consider a Banach subspace
$L^s(\Omega ,{\sf F},\lambda ;C^0_0(T,H))$
of $L^s(\Omega ,{\sf F},\lambda ;C^0(T,H))$,
where $\Omega \ne \emptyset $.
\par If $\mu $ is not a Gaussian measure on $C^0_0(T,H)$
and a stochastic process $w$ satisfies conditions $(i-iii)$, then
it is called the (non-Archimedean) stochastic process (see \S I.4.2).
If $\mu $ is the standard Gaussian measure on $C^0_0(T,H)$,
then the Wiener process is called standard
(see also Theorem 3.23, Lemmas 2.3, 2.5, 2.8 and
\S 3.30 in \cite{lu6}).
\par {\bf 3.3. Remark.} In Part I the non-Archimedean analogs 
of the It$\hat o$ formula were proved.  
In the particular case $H=\bf K$
we have $a\in L^s(\Omega ,{\sf F},\lambda ;C^0(T,{\bf K}))$,
$E\in L^r(\Omega ,{\sf F},\lambda ;C^0(T,{\bf K}))$,
$f\in C^n(T\times {\bf K},Y)$ and $w \in 
L^q(\Omega ,{\sf F},\lambda ;C^0_0(T,{\bf K}))$
are functions (see \S \S 4.2, 4.6 \cite{lust1} and \S 3.2), so that 
$${\hat P}_{u^{b+m-l},w(u,\omega )^l}[(\partial ^{m+b}f/
\partial u^b\partial x^m)(u,\xi (u,\omega ))
\circ (I^{\otimes b}\otimes a^{\otimes (m-l)}
\otimes E^{\otimes l})]|_{u=t}=$$ 
$$\sum_j (\partial ^{m+b}f/\partial u^b\partial x^m)(t_j,\xi (t_j,\omega ))
[t_{j+1}-t_j]^{b+m-l}a(t_j,\omega )^{k-l}
[E(t_j,\omega )(w(t_{j+1},\omega )-w(t_j,\omega ))]^l$$ for each $m+b\le n$,
where $t_j=\sigma _j(t)$,
$a(t,\omega ),$ $E(t,\omega )$ and $w(t,\omega )\in \bf K$, 
that is $a, E, w$ commute. In particular 
$\tilde P^m_{u,0}f(u)=\sum_{k=1}^m(k!)^{-1}{\hat P}_{u^k}f^{(k)}(u)$,
that is $\tilde P^m_{u,0}f(u)|_{u=t}=\tilde P_{m+1}f'(t),$ where
$\tilde P_{m+1}: C^m(T,{\bf K})\to C^{m+1}(T,{\bf K})$ 
is the Schikhof linear continuous antiderivation operator 
(see for comparison \S 80 \cite{sch1}).
\par In the classical case measures are real-valued and 
functions $\phi $ are with values in Banach spaces over
$\bf R$ or $\bf C$.
But in the considered here case measures are real-valued 
and functions are with values in Banach spaces
over non-Archimedean fields $\bf K$, so the mean value $M\| f \| $
is real and not with values in $\bf K$. This leads to differences 
with the classical case, in particular formula 
$M[(\int_S^T\phi (t,\omega )dB_t(\omega ))^2]=
M[\int_S^T\phi (t,\omega )^2dt]$ (see Lemma 3.5 \cite{oksen}) is not valid,
but there exists its another analog.
Let $X$ be a locally compact Hausdorff space and
$BC_c(X,H)$ denotes a subspace of $C^0(X,H)$ consisting
of bounded continuous functions $f$ such that for each
$\epsilon >0$ there exists a compact subset $V\subset X$
for which $\| f(u) \|_H <\epsilon $ for each $u\in X\setminus V$.
In particular for $X\subset \bf K$, $e^*\in H^*$ 
and a fixed $t\in X$ in accordance with Theorem 7.22 \cite{roo}
there exists a $\bf K$-valued tight measure 
$\mu _{t,\omega ,e^*,b,k}$ on the $\sigma $-algebra
$Bco(X)$ of clopen subsets in $X$ such that
$e^*{\hat P}_{u^b,w^k}\psi (u,x,\omega )\circ (I^{\otimes b}\otimes
E^{\otimes k})|_{u=t}=
\int_X\psi (u,E(u,\omega )w(u,\omega ),\omega )
\mu _{t,\omega ,e^*,b,k}(du)$ for each
$\psi \in L^r(\Omega ,{\sf F},\lambda ; BC_c(X,L_k(H^{\otimes k},H)))$ 
and $E\in L^q(\Omega ,{\sf F},\lambda ; BC_c(X,L(H)))$,
where $H^*$ is a topologically conjugate space, $1\le r, q \le \infty $,
$1/r+1/q\ge 1$.
\par If $\chi _{\gamma }: {\bf K}\to S^1:=\{ z\in {\bf C}: |z|=1 \}$
is a continuous character of $\bf K$ as the additive group,
then $M\chi _{\gamma }((e^*{\hat P}_{u^b,w^k}\psi (u,x,\omega )\circ 
(I^{\otimes b}\otimes E^{\otimes k})|_{u=t})^l)
=\prod_jM\chi _{\gamma }((e^*\psi (t_j,x,\omega )[t_{j+1}-t_j]^b\circ (
1^{\otimes b}\otimes 
(E(t_j,\omega )[w(t_{j+1},\omega )-w(t_j,\omega )])^{\otimes k})^l)$
due to Condition $I.4.2.(i).$ For $\psi $ independent from $x$,
$l=1$, $k=2$, $b=0$, $E=1$ and $H=\bf K$
(so that $e^*=1$) it takes a simpler form, which can be considered
as another analog of the classical formula. For the 
evaluation of appearing integrals tables from
\S 1.5.5 \cite{vla3} can be used.
Another important result is the following theorem.
\par {\bf Theorem.} {\it Let $\psi \in L^2(\Omega ,{\sf
F},\lambda ;C^0(T,L(H))),$  $w\in L^2(\Omega ,{\sf F},\lambda ;
C^0_0(T,H))$ be the stochastic process on the Banach space $H$ over $\bf K$.
Then there exists a function $\phi \in C^0(T,H)$ 
such that $M\chi _{\gamma }(g \hat P_{w(u,\omega )}\psi (u,\omega )\circ
I|_{u=t})=\hat \mu (\gamma g\hat P_u\phi (u)|_{u=t})$ for each
$\gamma \in \bf K$ and each $t\in T$ and for each $g\in H^*$.}
\par {\bf Proof.} Let $t\in T$ and $t_j=\sigma _j(t)$, where $\sigma _j$
is the approximation of the identity in $T$, $F_{a,b}(w):=
w(a,\omega )-w(b,\omega )$ for $a, b \in T$
(see \S I.2.1 \cite{lust1} and \S 3.2). In view of
Conditions $I.4.2.(i,ii)$ and the Hahn-Banach theorem
(see \cite{roo}) there exists a projection operator $Pr_g$
such that $\hat \mu ^{(F_{a,b}gE)}(h)=\hat \mu ^{(F_{a,b}Pr_g)}
(Pr_gEh)$, since $F_{a,b}ghEw=
ghE(w(a,\omega )-w(b,\omega ))=hgE F_{a,b}w$
for each $a, b \in T$ and for each $h\in \bf K$,
where $\hat \mu $ is the characteristic functional of the measure $\mu $
corresponding to $w$, that is,
$\hat \mu (g):=\int_{C^0_0(T,H)}\chi _g(y)\mu (dy)$, where
$g\in C^0_0(T,H)^*$, $\chi _g: C^0_0(T,H)\to \bf C$ is the character
of $C^0_0(T,H)$ as the additive group, $E\in L(H)$,
$y\in C^0_0(T,H)$, $\mu $ is the Borel measure on $C^0_0(T,H)$ (see also
\S I.3.6).
The random variable $E(w(a,\omega )-w(b,\omega ))$
has the distribution $\mu ^{F_{a,b}E}$ for each $a\ne b \in T$
and $E\in L(H)$.
On the other hand the projection
operator $Pr_e$ commutes with the antiderivation operator
$\hat P_u$ on $C^0(T,H)$, where $(Pr_ef)(t):=Pr_ef(t)$ is defined pointwise
for each $f\in C^0(T,H)$.
In $L^2(\Omega ,{\sf F},\lambda ;C^0(T,H))$ the family of step functions
$f(t,\omega )=\sum_{j=1}^nCh_{U_j}(\omega )f_j(t)$ is dense, where
$f_j\in C^0(T,H)$, $Ch_U$ is the characteristic function of $U\in \sf F$, 
$n\in \bf N$, since $\lambda (\Omega )=1$ and $\lambda $ is nonnegative.
 For each $t\in T$ there exists
$\lim_{j\to \infty }\psi (t_j,\omega ).(w(t_{j+1},\omega )-w(t_j,\omega ))$
in $L^2(\Omega ,{\sf F},\lambda ;H)$ (see Theorem I.2.14). 
\par If $A\in L(H)$, then 
\par $(i)$ $\chi _{\gamma }((g_1+g_2)Az)=\chi _{\gamma }(g_1Az)
\chi _{\gamma }(g_2Az)$ for each $g_1, g_2 \in H^*$ and $z\in H$,
\par $(ii)$ $\chi _{\gamma }(gA(z_1+z_2))=\chi _{\gamma }(gAz_1)
\chi _{\gamma }(gAz_2)$ for each $g\in H^*$ and $z_1, z_2 \in H$,
\par $(iii)$ $\chi _{\gamma }(agAz)=[\chi _{\gamma }(gAz)]^{\zeta (a)}$
for each $\{ (e,\gamma gAz) \} _p\ne 0$ and $a\in \bf K$,
where $\zeta (a):=\{ (e,\gamma agAz) \} _p/\{ (e,\gamma gAz) \} _p.$
On the other hand $A$ is completely defined by the family $ \{ e_i^*Ae_j: i,j
\in \alpha \} $, where $H=c_0(\alpha ,{\bf K})$, $e^*_i(e_j)=\delta _{i,j}$,
$e^*_i\in H^*$, $\{ e_j: j\in \alpha \} $ is the standard orthonormal base
of $H$. Hence the family $ \{ \chi _{\gamma }(ae^*_iAe_j): i, j \in \alpha ;
a\in {\bf K} \} $ completely characterize $A\in L(H)$ due to Equations
$(i-iii)$, when $\gamma \ne 0$.
\par For each $y\in H$ and each $\gamma \in \bf K$ the function
$M\chi _{\gamma }(g\psi (t,\omega )y)$ is continuous by $t\in T$, consequently,
there exists a continuous function $\phi : T\to H$ such that
$M\chi _{\gamma }(g\psi (t,\omega )y)=\chi _{\gamma }(g\phi (t)y)$
for each $y\in H$ and $t\in T$, since characters $\chi _{\gamma }$
are continuous from $\bf K$ to $\bf C$ and $\chi _{\gamma }(h)=\chi _1
(\gamma h)$ for each $0\ne \gamma \in \bf K$ and $h\in \bf K$ and the 
$\bf C$-linear span of the family $\{ \chi _{\gamma }: \gamma \in {\bf K} \} $
of characters is dense in $C^0({\bf K},{\bf C})$ 
by the Stone-Weierstrass theorem 
\cite{fell}. On the other hand, $\lim_{j\to \infty }
\chi _{\gamma }(\sum_{i=0}^ja_j)=\prod _{i=1}^{\infty }\chi _{\gamma }(a_i)$,
when $\lim_ja_j=0$ for a sequence $a_j$ in $\bf K$. 
Therefore, 
$$M\chi _{\gamma }(g\sum_{j=0}^{\infty }
\psi (t_j,\omega ).[w(t_{j+1},\omega )-w(t_j,\omega )])=
\prod_{j=0}^{\infty }\hat \mu (\gamma g\phi (t_j)(t_{j+1}-t_j))$$ 
$$=\hat \mu (\gamma g\hat P_u\phi (u)|_{u=t}) \quad \mbox{ for each }t\in T
\mbox{ and each }g\in H^*.$$
From the equality $\chi _{a+b}(c)=\chi _a(c)\chi _b(c)$
for each $a, b$ and $c\in \bf K$ the statement of this theorem follows for 
each $\gamma \in \bf K$.
\par {\bf 3.4.} {\bf Theorem.} {\it 
Let $a\in L^q(\Omega ,{\sf F},\lambda ;C^0(B_R,L^q(
\Omega ,{\sf F},\lambda ;C^0(B_R,H))))$
and $E\in L^q(\Omega ,{\sf F},\lambda ;C^0(B_R,L(L^q(
\Omega ,{\sf F},\lambda ;C^0(B_R,H))))),$
$a=a(t,\omega ,\xi )$, $E=E(t,\omega ,\xi )$, 
$t\in B_R,$ $\omega \in \Omega ,$
$\xi \in L^q(\Omega ,{\sf F},\lambda ;C^0(B_R,H))$
and $\xi _0 \in L^q(\Omega ,{\sf F},\lambda ;H),$
and $w \in L^q(\Omega ,{\sf F},\lambda ;C^0_0(B_R,H)),$
where $a$ and $E$ satisfy the local Lipschitz condition:
\par $(LLC)$ for each $0<r<\infty $ there exists $K_r>0$ such that
$\max ( \| a(t,\omega ,x)-a(t,\omega ,y) \|, 
\| E(t,\omega ,x)-E(t,\omega ,y) \| )\le K_r \| x-y \| $
for each $x, y \in B(C^0(B_R,H),0,r)$ and $t\in B_R$, 
$\omega \in \Omega $, $1\le q \le \infty $.
Then the stochastic process of the following type:
\par $(i)$ $\xi (t,\omega )=\xi _0(\omega )+
(\hat P_ua)(u,\omega ,\xi )|_{u=t}+({\hat P}_{w(u,\omega )}E)(
u,\omega ,\xi )|_{u=t}$ has the unique solution.}
\par {\bf Proof.} We have $\max ( \| a(x)-a(y) \|^g, \| E(x)-E(y) \|^g )
\le K\| x-y \|^g $, hence $\max ( \| a(x) \|^g, \| E(x) \|^g )\le K_1
(\| x \|^g +1)$ for each $x, y\in H$ and for each $1\le g<\infty $
and each $t\in B_R$ and each $\omega \in \Omega $, 
where $K$ and $K_1$ are positive 
constants, $a(x)$ and $E(x)$ are short notations of $a(t,\omega ,x)$ and
$E(t,\omega ,x)$ for $x=\xi (t,\omega )$
respectively. For solving equation $(i)$ we use iterations:
\par $X_0(t)=x$,..., $X_n(t)=x+\hat P_ua(X_{n-1}(u))|_{u=t}+
{\hat P}_wE(X_{n-1}(u))|_{u=t}$, consequently,
$X_{n+1}-X_n(t)=I_1(t)+I_2(t)$, where $I_1(t)=\hat P_u[a(X_n(u))-a(X_{n-1}(u))]
|_{u=t}$, $I_2(t)={\hat P}_w[E(X_n(u))-E(X_{n-1}(u))]|_{u=t}$,
$x(t)$ and $X_n(t)$ are short notations of $x(t,\omega )$ 
and $X_n(t,\omega )$ respectively.
Let $M\eta $ be a mean value of a real-valued distribution
$\eta (\omega )$ by $\omega \in \Omega $, where $(\Omega ,
{\sf F},\lambda )$ is the probabilty space, then
$M \| P_u[a(X_n(u))-a(X_{n-1}(u))]|_{u=t} \|^g \le K(M \| P_t \|^g)
M\sup_u\| X_n(u)-X_{n-1}(u)\| ^g$, where 
$X_n\in L^q(\Omega ,{\sf F},\lambda ;C^0_0(B_R,H))$
for each $n$, since $|\lambda |(\Omega )=1$ and $\| x \| _{\infty }
=\sup_{1\le g<\infty }\| x \| _g=ess-\sup_{\omega \in \Omega } 
\| x(\omega ) \| _H$ for $x\in L^{\infty }(\Omega ,{\sf F},\lambda ;H)$. 
While $1\le q<\infty $ we put $g=q$, for $q=\infty $ we take
$ess-\sup $. Also 
$$M\| {\hat P}_w[E(X_n(u))-E(X_{n-1}(u))]|_{u=t} \| ^g
\le K \| {\hat P}_w\| ^g M\sup_u\| X_n(u)-X_{n-1}(u)\| ^g$$
$$\le (K\| {\hat P}_w \| ^g)^lM\sup_u \| X_{n-l+1}(u)-X_{n-l}(u) \| ^g$$
in particular for $l=n-1$. On the other hand,
$$X_1(t)=x(t)+\hat P_ua(x(u))|_{u=t}+{\hat P}_wE(x(u))|_{u=t},$$
consequently, $\| X_1(t)-X_0(t) \| ^g\le \max (\| [P_ua(x(u))]|_{u=t} \| ^g,
\| {\hat P}_wE(x(u))|_{u=t_0} \| ^g)$, where $w(0)=0$, $\hat P_ua(u)|_{u=t_0}=0$,
${\hat P}_wE|_{u=t_0}=0$. For each $\epsilon >0$ there exists
$B_{\epsilon }\subset B_R$ such that $K \| {\hat P}_w|_{B_{\epsilon }}
\| ^g<1$ and $K \| P_t|_{B_{\epsilon }} \| ^g <1$.
Therefore, there exists the unique solution on each
$B_{\epsilon }$, since $\sup_u\| X_1(u)-X_0(u) \| ^g<\infty $
and $\lim_{l\to \infty }(K \| {\hat P}_w|_{B_{\epsilon }} \| ^g)^lC=0$,
$\lim_{l\to \infty }(K \| P_t|_{B_{\epsilon }} \| ^g )^lC=0$,
hence there exists $\lim_{n \to \infty }X_n(t)=X(t)=\xi (t,\omega )
|_{B_{\epsilon }}$, where
$C:=M\sup_{u\in B_{\epsilon }}\| X_1(u)-X_0(u) \| ^g
\le \max (\| {\hat P}_w \| ^g, \| P_t \| ^g)(\| x \|_{C^0}^g+1)K<\infty ,$ 
here $B_{\epsilon }$ is an arbitrary ball of radius 
$\epsilon $ in $B_R$, $t\in B_{\epsilon }$. 
\par If $X^1$ and $X^2$ are two solutions, then
$X^1-X^2=:\psi =\sum_{j=1}^nC_jCh_{B({\bf K},x_j,r_j)}$, where
$n\in \bf N$, $C_j\in \bf K$, $T=B_R$, since $B_R$ has 
a disjoint covering by balls $B({\bf K},x_j,r_j)$,
on each such ball there exists the unique solution
with a given initial condition on it (that is, in a 
chosen point $x_j$ such that $C_j$ and 
$B({\bf K},x_j,r_j)$ are independent from $\omega $). Therefore,
$\psi =\hat P_u[a(u,X^2)-a(u,X^1)]|_{u=t}+{\hat P}_w[E(u,X^2)-E(u,X^1)]|_{u=t}$,
hence $\Phi ^1\psi (t_i;1;t_{i+1}-t_i)=
[a(t_i,X^2(t_i))-a(t_i,X^1(t_i))+[E(t_i,X^2(t_i))-E(t_i,X^1(t_i))][w(t_{i+1})
-w(t_i)]/(t_{i+1}-t_i)$ for each $t_i\ne t_{i+1}$, $t_i=\sigma _i(t)$
due to Condition $I.2.1.(ii)$, where $w(t)$ is the short notation of
$w(t,\omega )$.
The term $(\Phi ^1w)(t_i;1;t_{i+1}-t_i)=[w(t_{i+1})-w(t_i)]/(t_{i+1}-t_i)$
has the infinite-dimensional over $\bf K$
range in $C^0(B_R^2\setminus \Delta ,H)$
for each $\omega \in \Omega $, where $\Delta := \{ (u,u): u\in B_R \} $.
If $(\Phi ^1w)=0$, then $a(t,X^2(t,\omega ))-a(t,X^1(t,\omega ))=0$. 
If $a(t,X^2)=a(t,X^1)$ for each $t$ and almost all $\omega $, then
${\hat P}_w[E(t,X^2(t))-E(t,X^1(t))]=0$ which is possible 
only for $\psi =0$.
If $a(t,X^2)\ne a(t,X^1)$ and the function $\psi $ is locally constant
by $t$ and independent from $\omega $, then 
$\hat P_u[a(u,X^1+g)-a(u,X^1)]|_{u=t}+{\hat P}_w[E(u,X^1+g)-E(u,X^1)]|_{u=t}$
is locally constant by $t$ and independent from $\omega $ only for 
$[a(u,X^2)-a(u,X^1)]=0$ and $[E(u,X^2)-E(u,X^1)]=0$ due to definitions 
of $\hat P_u$ and ${\hat P}_w$, hence $\psi =0$, since it is evident
for $a(u,X)$ and $E(u,X)$ depending on $X$ locally polynomially 
or polyhomogeneously for each $u$, but such locally
polynomial or polyhomogeneous functions by $X$ are dense in 
$$L^q(\Omega ,{\sf F},\lambda ;C^0(B_R,L^q(
\Omega ,{\sf F},\lambda ;C^0(B_R,H))))\mbox{ and}$$ 
$$L^q(\Omega ,{\sf F},\lambda ;C^0(B_R,L(L^q(
\Omega ,{\sf F},\lambda ;C^0(B_R,H)))))\mbox{ respectively}.$$
\par {\bf 3.5. Theorem.} {\it Let $a\in 
L^q(\Omega ,{\sf F},\lambda ;C^0(B_R,L^q(
\Omega ,{\sf F},\lambda ;C^0(B_R,H))))$
and $E\in L^q(\Omega ,{\sf F},\lambda ;C^0(B_R,L(L^q(
\Omega ,{\sf F},\lambda ;C^0(B_R,H))))),$
$a=a(t,\omega ,\xi )$, $E=E(t,\omega ,\xi )$, 
$t\in B_R,$ $\omega \in \Omega ,$
$\xi \in L^q(\Omega ,{\sf F},\lambda ;C^0(B_R,H))$
and $\xi _0 \in L^q(\Omega ,{\sf F},\lambda ;H),$
where $a$ and $E$ satisfy the local Lipschitz condition
(see 3.4.(LLC)).
A stochastic process of the type
\par $(i)$ $\xi (t,\omega )=\xi _0(\omega )+$
$\sum_{m+b=1}^{\infty }\sum_{l=0}^m(
{\hat P}_{u^{b+m-l},w(u,\omega )^l}[a_{m-l+b,l}
(u,\xi (u,\omega ))\circ (I^{\otimes b}\otimes 
a^{\otimes (m-l)}\otimes E^{\otimes l})])
|_{u=t}$  \\
such that
$a_{m-l,l} \in C^0(B_{R_1}\times B(L^q(\Omega ,{\sf F},\lambda ;C^0(B_R,H)),
0,R_2),L_m(H^{\otimes m};H))$ (continuous and bounded on its domain)
for each $n, l,$ $0<R_2<\infty $ and
\par $(ii)$ $\lim_{n\to \infty } \sup_{0\le l\le n}\|a_{n-l,l}
\|_{C^0(B_{R_1}\times B(L^q(\Omega ,{\sf F},\lambda ;C^0(B_R,H)),
0,R_2),L_n(H^{\otimes n},H))}=0$ 
for each $0<R_1\le R$ when $0<R<\infty $, or each $0<R_1<R$
when $R=\infty $, for each $0<R_2<\infty .$
\par Then $(i)$ has the unique solution in $B_R$.}
\par {\bf Proof.} Let $X_0(t)=x$,..., 
\par $$X_n(t)=x+
\sum_{m+b=1}^{\infty }\sum_{l=0}^m(
{\hat P}_{u^{b+m-l},w(u,\omega )^l}[a_{m-l+b,l}
(u,X_{n-1}(u,\omega ))\circ (I^{\otimes b}\otimes a^{\otimes (m-l)}
\otimes E^{\otimes l})])
|_{u=t},$$ consequently,
$$X_{n+1}-X_n(t)=
\sum_{m+b=1}^{\infty }\sum_{l=0}^m(
{\hat P}_{u^{b+m-l},w(u,\omega )^l}[a_{m-l+b,l}
(u,X_n(u))-a_{m-l+b,l}(u,X_{n-1}(u))]$$
$\circ (I^{\otimes b}\otimes a^{\otimes (m-l)}\otimes E^{\otimes l})])
|_{u=t},$ \\
where in general ${\hat P}_{a(u,\xi )}1|_{u=t}=
a(t,\xi (t,\omega ))-a(t_0,\xi (t_0,\omega ))
\ne \hat P_ua(u,\xi )=\sum_ja(t_j,\xi (t_j,\omega ))[t_{j+1}-t_j]$,
$t_j=\sigma _j(t)$ for each $j=0,1,2,...$. Then
$$M \| {\hat P}_{u^{b+m-l},w(u,\omega )^l}[a_{m-l+b,l}
(u,X_n(u))-a_{m-l+b,l}(u,X_{n-1}(u))]|_{(B_{R_1}\times B(L^q,0,R_2))}
\circ ( $$ 
$$I^{\otimes b}\otimes a^{\otimes (m-l)}\otimes E^{\otimes l})])
|_{u=t} \| ^g \le K(M \| {\hat P}_{u^{b+m-l},w(u,\omega )^l} \| ^g)
\| a_{m-l+b,l} |_{(B_{R_1}\times B(L^q,0,R_2))}\| ^g$$ 
$$ (M \sup_u \| X_n(u)-X_{n-1}(u) \| ^g)
(M\sup_u \| a \| ^{m-l})(M\sup_u \| E \| ^l),$$
where $X_n\in C^0_0(B_R,H)$
for each $n$, $K$ is the same constant as in \S 3.4, $1\le g<\infty $. 
On the other hand,
$$X_1(t)=x(t)+
\sum_{m+b=1}^{\infty }\sum_{l=0}^m(
{\hat P}_{u^{b+m-l},w(u,\omega )^l}[a_{m-l+b,l}
(u,x(u))\circ (I^{\otimes b}\otimes a^{\otimes (m-l)}
\otimes E^{\otimes l})])|_{u=t},$$ consequently,
$$\| X_1(t)-X_0(t) \| ^g \le \sup_{m,l,b} (\| 
{\hat P}_{u^{b+m-l},w(u,\omega )^l}[a_{m-l+b,l}
(u,x(u))\circ (I^{\otimes b}\otimes a^{\otimes (m-l)}
\otimes E^{\otimes l})])|_{u=t} \| ^g .$$
Due to Condition $(ii)$ 
for each $\epsilon >0$ and $0<R_2<\infty $ there exists
$B_{\epsilon }\subset B_R$ such that 
$$K \sup_{m,l,b} (\| 
{\hat P}_{u^{b+m-l},w(u,\omega )^l}|_{B_{\epsilon }} 
[a_{m-l+b,l}(u,*)|_{(B_{\epsilon }\times B(L^q,0,R_2))}\circ (I^{\otimes b}\otimes a^{\otimes (m-l)}
\otimes E^{\otimes l})]) \| ^g =:c<1.$$
Therefore, there exists the unique solution on each
$B_{\epsilon }$, since $\sup_u\| X_1(u)-X_0(u) \| <\infty $
and $\lim_{l\to \infty }c^lC=0$ for each $C>0$,
hence there exists $\lim_{n \to \infty }X_n(t)=X(t)=\xi (t,\omega )
|_{B_{\epsilon }}$, where
$C:=M\sup_{u\in B_{\epsilon }}\| X_1(u)-X_0(u) \| ^g
\le (c+1)K<\infty ,$ 
here $B_{\epsilon }$ is an arbitrary ball of radius 
$\epsilon $ in $B_R$, $t\in B_{\epsilon }$. 
\par If $X^1$ and $X^2$ are two solutions, then
$X^1-X^2=:\psi =\sum_{j=1}^nC_jCh_{B({\bf K},x_j,r_j)}$ as in \S 3.4.
If $S$ is a polyhomogeneous function, then there exists
$n=deg(S)<\infty $ such that differentials $D^mS=0$ for each $m>n$,
but its antiderivative $\hat P$ has $D^{n+1}{\hat P}S\ne 0$.
If $\| S_1 \| > \| S_2 \| $, then $\| {\hat P}S_1 \| >
\| {\hat P}S_2 \| $, which we can apply to a convergent series
considering terms $\| D^m{\hat P} S \| (mod \mbox{ }p^k)$
for each $k\in \bf N$. Therefore, 
\par $\psi =
\sum_{m+b=1}^{\infty }\sum_{l=0}^m(
{\hat P}_{u^{b+m-l},w(u,\omega )^l}[a_{m-l+b,l}
(u,X^2)-a_{m-l+b,l}(u,X^1)]
\circ (I^{\otimes b}\otimes a^{\otimes (m-l)}
\otimes E^{\otimes l})])|_{u=t}$, 
where the function $\psi $ is locally constant
by $t$ and independent from $\omega $, 
hence $\psi =0$, since it is evident
for $a(u,X)$ and $E(u,X)$ and $a_{k-l,l}(u,X)$ depending on $X$ 
locally polynomially 
or polyhomogeneously for each $u$, but such locally
polynomial or polyhomogeneous functions by $X$ are dense in 
$$L^q(\Omega ,{\sf F},\lambda ;C^0(B_R,L^q(
\Omega ,{\sf F},\lambda ;C^0(B_R,H))))\mbox{ and}$$ 
$$L^q(\Omega ,{\sf F},\lambda ;C^0(B_R,L(L^q(
\Omega ,{\sf F},\lambda ;C^0(B_R,H)))))\mbox{ and}$$ 
$$C^0(B_{R_1}\times B(L^q(\Omega ,{\sf F},\lambda ;C^0(B_R,H)),0,R_2),
L_k(H^{\otimes k};H))$$ 
respectively.
\par {\bf 3.6. Proposition.} {\it Let $\xi $ be the Wiener
process given by Equation $3.4.(i)$ with the $1$-Gaussian measure
associated with the operator $\tilde P^1\tilde P^0$ as in \S 2.4
and let also $\max (\| a(t,\omega ,x)-a(v,\omega ,x) \| , 
\| E(t,\omega ,x)-E(v,\omega ,x) \| ) \le |t-v|(C_1+C_2 \| x \| ^b)$
for each $t$ and $v\in B({\bf K},t_0,R)$ $\lambda $-almost everywhere
by $\omega \in \Omega $,
where $b$, $C_1$ and $C_2$ are non-negative constants.
Then $\xi $ with probability $1$ has a $C^2$-modification
and  $q(t)\le \max \{ \| \xi _0 \| ^s,
|t-t_0|(C_1+C_2q(t)) \} $
for each $t\in B({\bf K},t_0,R)$, where
$q(t):= \sup_{|u-t_0|\le |t-t_0|}M \| \xi (t,\omega )\| ^s $
and ${\bf N}\ni s\ge b\ge 0$.}
\par {\bf Proof.} For the following function
$f(t,x)=x^s$ in accordance with Theorem I.4.6 \cite{lust1} we have
$f(t, \xi (t,\omega ))=f(t_0,\xi _0)+$
$$+\sum_{k=1}^s\sum_{l=0}^k{k\choose l}(
\hat P_{u^{k-l},w(u,\omega )^l}[({s\choose k}\xi (t,\omega )^{s-k}
(u,\xi (u,\omega ))\circ (a^{\otimes (k-l)}\otimes E^{\otimes l})])
|_{u=t},$$ hence 
$M \| \xi (t,\omega ) \| ^s\le 
\max ( \| \xi _0 \| ^s,  |t-t_0| 
d(\hat P^s_*) (C_1+C_2 \sup_{|u-t_0|\le |t-t_0|}
M \| \xi (u,\omega )\| ^s)$, since $|t_j-t_0|
\le |t-t_0|$ for each $j\in \bf N$ and
$M \| \xi (t,\omega ) -\xi (v,\omega )\| ^s\le 
|t-v| (1+C_1+C_2 d(\hat P^s_*) \sup_{|u-t_0|\le \max( |t-t_0|, |v-t_0|)}
M \| \xi (u,\omega )\| ^s)$, since $|t_j-v_j|\le |t-v|+\rho ^j$
for each $j\in \bf N$, where $0<\rho <1$, 
$$d( \hat P^s_* ) :=\sup_{a\ne 0, E\ne 0, f \ne 0}
\max_{s\ge k\ge l\ge 0}
\| (k!)^{-1}{k\choose l}{\hat P}_{u^{k-l},w^l}(\partial ^kf
/\partial ^kx ) \circ (a^{\otimes (k-l)}\otimes E^{\otimes l}) \| /$$ 
$$(\| a \|_{C^0(B_R,H)}^{k-l} \| E \|_{C^0(B_R,L(H))}^l 
\| f \|_{C^s(B_R,H)} ),$$ 
hence $d( \hat P^s_* )\le 1$, since
$f\in C^s$ as a function by $x$ and $(\bar \Phi ^sg)(x;h_1,...,h_s;0,...,0)=
D^s_xg(x).(h_1,...,h_s)/s!$ for each $g\in C^s$ and due to 
the definition of $\| g \| _{C^s}$. Considering in particular 
polyhomogeneous $g$ on which $d(\hat P^s_*)$ takes its maximum value 
we get $d(\hat P^s_*)=1$. Since $P(C^2)=1$ for the Markov measure
$P$ induced by the transition measures 
$P(v,x,t,S):=\mu ^{F_{t,v}}(S|\xi (v)=x)$ for $t\ne v$ 
of the non-Archimedean Wiener process (see \S 2.2), then
$\xi $ has with the probability $1$ a $C^2$-modification. 
\par {\bf Note.} If to consider a general stochastic process
as in \S I.4.3, then from the proof of Proposition 3.6 it follows, that
$\xi $ with the probability $1$ has a modification in the space
$J(C^0_0(T,H))$, where $J$ is a nondegenerate correlation operator
of the product measure $\mu $ on $C^0_0(T,H)$.
\par {\bf 3.7. Proposition.} {\it Let $\xi $ be a stochastic
process given by Equation $3.4.(i)$ and $\max 
(\| a (t,\omega ,x_1) -a(v,\omega ,x_2) \| , \| E (t,\omega ,x_1) 
-E(v,\omega ,x_2) \| ) 
\le |t-v| (C_1+C_2 \| x_1-x_2 \| ^b)$ for each $t$ and $v\in 
B({\bf K},t_0,R)$ $\lambda $-almost everywhere by $\omega \in \Omega $,
where $b$, $C_1$ and $C_2$ are non-negative constants.
Then two solutions $\xi _1$ and $\xi _2$ with initial conditions 
$\xi _{1,0}$ and $\xi _{2,0}$ satisfy the following inequality:
$y(t)\le \max \{ \| \xi _{1,0}-\xi _{2,0} \| ^s,
|t-t_0|(C_1+C_2y(t)) \} $
for each $t\in B({\bf K},t_0,R)$, where
$y(t):= \sup_{|u-t_0|\le |t-t_0|}M \| \xi _1(t,\omega )- 
\xi _2(t,\omega )\| ^s $
and ${\bf N}\ni s\ge b\ge 0$.}
\par {\bf Proof.} From \S 3.6 it follows, that
$M \| \xi _1(t,\omega ) - \xi _2(t,\omega )\| ^s\le 
|t-t_0| (C_1+C_2 \sup_{|u-t_0|\le |t-t_0|}
M \| \xi _1(u,\omega )-\xi _2(u,\omega )\| ^s)$, since $d(\hat P^s_*)\le 1$.
\par {\bf 3.8. Remark.} Let $X_t=X_0+\hat P_ta+{\hat P}_wv$
and $Y_t=Y_0+\hat P_tq+{\hat P}_ws$ be two stochastic processes
corresponding to $E=I$ and a Banach algebra $H$ over $\bf K$
in \S I.4.6 \cite{lust1,boua,roo}. Then $X_uY_u-X_tY_t=
(X_u-X_t)(Y_u-Y_t)+X_t(Y_u-Y_t)+(X_u-X_t)Y_t$, where $u, t\in T$.
Hence $d(X_tY_t)=X_tdY_t+(dX_t)Y_t+(dX_t)(dY_t)$.
Therefore, 
$${\hat P}_{X_t}Y_t=X_tY_t-X_0Y_0-
{\hat P}_{Y_t}X_t-{\hat P}_{(X_t,Y_t)}1,$$
which is the non-Archimedean analog of the integration by parts formula,
where in all terms $X_t$ is displayed on the left from $Y_t$.
For two $C^1$ functions $f$ and $g$ we have $(fg)'=f'g+fg'$ or
$d(fg)=gdf+fdg$, that is terms with $(dt)(dt)$ are absent,
consequently, $(dt)(dt)=0$. In a particular case $X_t=Y_t=w_t$
this leads two $w_t^2-w_0^2-2{\hat P}_{w_t}w_t={\hat P}_{(w_t,w_t)}
1$, where the last term corresponds two $(dw_t)(dw_t)\ne 0$.
This means that 
$$d(w^2)=2wdw+(dw)(dw).$$ 
For $X_t=w_t$ and $Y_t=t$ 
the integration by parts formula gives ${\hat P}_{w_t}t=
w_tt-\hat P_tw_t-{\hat P}_{(t,w_t)}1$. Such that
${\hat P}_{(t,w_t)}1=\sum_jt_j[w_{t_{j+1}}-w_{t_j}]-w_tt+
\sum_jw_{t_j}[t_{j+1}-t_j]\ne 0$, for example, for $t=1$, $w\in C^0_0(T,H)$,
$T=\bf Z_p$ and $t_0=0$
this gives ${\hat P}_{(t,w_t)}1=w_1-w_0=w_1$. Therefore,
$(dt)(dw_t)\ne 0$, that is the important difference of the non-Archimedean 
and classical cases (see for comparison Exer. 4.3 and Theorem 
4.5 \cite{oksen}).
\par If $H$ is a Banach space over the local field $\bf K$
and $f(x,y)=x^*y$ is a $\bf K$-bilinear 
functional on it, where $x^*$ is an image of $x\in H$ 
under an embedding $H\hookrightarrow H^*$ associated with the standard 
orthonormal base $\{ e_j \} $ in $H$, then 
$${\hat P}_{X^*_t}Y_t=X^*_tY_t-X^*_0Y_0-
{\hat P}_{Y^*_t}X_t-{\hat P}_{(X^*_t,Y_t)}1,$$
hence $d(X^*_tY_t)=X^*_tdY_t+(dX^*_t)Y_t+(dX^*_t)(dY_t)$
and $d(w^*w)=w^*dw+(dw^*)w+(dw^*)(dw).$
\par {\bf 3.9. Definition.} If $\xi (t,\omega )
\in L^q(\Omega ,{\sf F},\lambda ;C^0(B_R,H))=:Z$
is a stochastic process and $T(t,s)$ is a 
family of bounded linear operators satisfying the following Conditions $(i-iv):$ 
\par $(i)$ $T(t,s): H_s\to H_t$, where $H_s:=
L^q(\Omega ,{\sf F},\lambda ;C^0(B({\bf K},0,|s|),H)),$
\par $(ii)$ $T(t,t)=I$, 
\par $(iii)$ $T(t,s)T(s,v)=
T(t,v)$ for each $t, s, v\in B_R$,
\par $(iv)$ $M_s \{ \| T(t,s) \eta \| ^q_H \}
\le C \| \eta \|^q_H$ for each $\eta \in H_s$,
where $C$ is a positive nonrandom constant,
$1\le q \le \infty $, then $T(t,s)$ 
is called a multiplicative operator functional
of the stochastic process $\xi $. 
\par If $T(t,s;\omega )$ is a system of random variables
on $\Omega $ with values in $L(H)$, satisfying almost surely
Conditions $(i-iii)$ and uniformly by $t, s \in B_R$
Condition $(iv)$ such that
\par $(v)$ $(T(t,s)\eta )(\omega )=T(t,s;\omega )\eta (\omega )$,
then such multiplicative operator functional is called homogeneous.
An operator
\par $(vi)$ $A(t)=\lim_{s\to 0}[T(t,t+s)-I]/s$
is called the generating operator of the evolution family $T(t,v)$.
If $T(t,v)=T(t,v;\omega )$ depends on $\omega $, then $A(t)=A(t;\omega )$
is also considered as the random variable on $\Omega $ (depending
on the parameter $\omega $) with values in $L(H)$.
\par {\bf 3.10. Remark.} Let $A(t)$ be a linear continuous operator
on a Banach space $Y$ over $\bf K$ such that it depends
strongly continuously on $t\in B({\bf K},0,R)$, that is 
$A(t)y$ is continuous by $t$ for each chosen $y\in Y$
and $A(t)\in L(Y)$. Then the solution of the differential equation
\par $(1)\mbox{ }dx(t)/dt=A(t)x(t)\mbox{, }x(s)=x_0\mbox{ has a solution }$
\par $(2)\mbox{ }x(t)=U(t,s)x(s)\mbox{, where }U(t,s)\mbox{ is
a generating operator such that}$
\par $(3)\mbox{ }U(t,s)=I+{\hat P}_uA(u)U(u,s)|^{u=t}_{u=s},$
though $x(t)$ may be non-unique, where $x(s)=x_0$ is an initial 
condition, $x, t \in B({\bf K},0,R)$. 
The solution of Equation $(3)$ exists using the method
of iterations (see \S 3.4).
Indeed, in view of Lemma I.2.3 \cite{lust1} $U(s,s)=I$ and
\par $(4)$ $dx(t)/dt
=\partial U(t,s)x(s)/\partial t=A(t)U(t,s)x(s)=
A(t)x(t)$. If to consider a solution of the antiderivational equation
\par $(5)\mbox{ }V(t,s)=I+{\hat P}_uV(t,u)A(u)|^{u=t}_{u=s}$, then it is
a solution of the Cauchy problem
\par $(6)$ $\partial V(t,s)/\partial s=-V(t,s)A(s)$, $V(t,t)=I$.
Therefore, $\partial [V(t,s)U(s,v)]/\partial s
=-V(t,s)A(s)U(s,v)+V(t,s)A(s)U(s,v)=0,$
hence $V(t,s)U(s,v)$ is not dependent from $s$, consequently,
there exist $U$ and $V$ such that
\par $(7)$ $V(t,s)=U(t,s)$ for each $t, s \in B({\bf K},0,R)$. 
From this it follows, that 
\par $(8)$ $U(t,s)U(s,u)=U(t,u)$ for each $s, u, t \in B({\bf K},0,R)$.
In particular, if $A(t)=A$ is a constant operator, then
there exists a solution $U(t,s)=EXP((t-s)A)$ (see about $EXP$
in Proposition 45.6 \cite{sch1}).
Equation $(3)$ has a solution under milder conditions, for example, 
$A(t)$ is weakly continuous, that is $e^* A(t) \eta $ is 
continuous for each $e^*\in Y^*$ and $\eta \in Y$, then $e^*U(t,s)\eta $ is
differentiable by $t$ and $U(t,s)$ satisfies Equation $(4)$ in the weak 
sense and there exists a weak solution of $(5)$ coinciding with
$U(t,s)$. If to substitute $A(t)$ on another operator ${\tilde A}(t)$, then
for the corresponding evolution operator ${\tilde U}(t,s)$
there is the following inequality:
\par $(9)$ $ \| {\tilde U}(t,s)-U(t,s) \| \le M{\tilde M} 
\sup_{u\in B({\bf K},0,R)}\| {\tilde A}(u)-A(u) \| R$,
where $M:=1+\sup_{s, t \in B({\bf K},0,R)} \| U(t,s) \| $
and $\tilde M$ is for $\tilde U$.
\par {\bf Proposition.} {\it Let $B(t)$ and two sequences
$A_n(t)$ and $B_n(t)$ be given of 
strongly continuous on $B({\bf K},0,R)$ bounded linear operators
and ${\tilde U}(t,s)$ be evolution operators corresponding to 
${\tilde A}_n(t)=A_n(t)+B_n(t)$, where \\
$\sup_{n\in {\bf N}, 
u\in B({\bf K},0,R)} \| B_n(u) \| \le \sup_{u\in B({\bf K},0,R)}
\| B(u) \| =C< \infty $. If $MCR<1$, 
then there exists a sequence ${\tilde U}_n(t,s)$ which 
is also uniformly bounded. If there exists $U_n(t,s)$ strongly and 
uniformly converging to $U(t,s)$
in $B({\bf K},0,R)$, then ${\tilde U}_n(t,s)$ also can be chosen
strongly and uniformly convergent.}
\par {\bf Proof.} From the use of Equations $(3,8)$ iteratively for 
$U_n(\sigma _{j+1}(t),\sigma _j(t))$ and 
$U_n(\sigma _j(t),s)$ and also for $\tilde U_n$ and
taking $\tilde U_n-U_n$ it follows, that
\par $(10)$ ${\tilde U}_n(t,s)=
U_n(t,s)+{\hat P}_vU_n(t,v)B_n(v){\tilde U}_n(
v,s)|^{v=t}_{v=s}$ for each $n\in \bf N$. Therefore,
$\| {\tilde U}_n(t,s) \| \le M+MC \sup_v\| {\tilde U}_n(v,s) \| R$, hence
$\| {\tilde U}_n(t,s) \| \le M/[1-MCR]$, since $MCR<1$. 
If $\lim_nx_n=x$ in $Y$ and $U_n(t,s)x$ is uniformly convergent to
$U(t,s)x$, then for each $\epsilon >0$ there exist $\delta >0$ and 
$m\in \bf N$ such that $\sup_{t, s\in B({\bf K},0,R)}
\| U_n(t+h,s+v)x_n-U_n(t,s)x_n \| <\epsilon $
for each $n>m$ and $\max(|h|,|v|)<\delta $ due to Equality $(10)$.
\par {\bf 3.11. Proposition.} {\it Let $a$, $a_{m-l+b,l}$ and $E$ be the same 
as in \S 3.5. Then Equation $3.5.(i)$ has the unique solution $\xi $ in $B_R$
for each initial value $\xi (t_0,\omega )\in L^q(\Omega ,{\sf F},\lambda ;H)$
and it can be represented in the following form:
\par $(2)$ $\xi (t,\omega )=T(t,t_0;\omega )\xi (t_0;\omega )$,
where $T(t,v;\omega )$ is the multiplicative operator functional.}
\par {\bf Proof.} In view of Theorem 3.5, Definition 3.9, 
Remark and Proposition 3.10
with the use of a parameter $\omega \in \Omega $
the statement of Proposition 3.11 follows.
\par {\bf 3.12.} Let now consider the case $J(C^0_0(T,H))\subset C^1(T,H)$ 
(see \S 3.6), for example, the standard Wiener process.
\par {\bf Corollary.} {\it Let a function $f(t,x)$ satisfies 
conditions of \S I.4.8 \cite{lust1}, then a generating operator of an evolution 
family $T(t,v)$ of a stochastic process 
$\eta =f(t,\xi (t,\omega ))$ is given by the following equation:
\par $(1)$ $A(t)\eta (t)=f'_t(t,\xi (t,\omega ))$
$+f'_x(t,\xi (t,\omega ))\circ a(t,\omega )+$
$$f'_x(t,\xi (t,\omega ))\circ E(t,\omega ){w'}_t(t,\omega )+
\sum_{m+b\ge 2, 0\le m\in {\bf Z}, 0\le b\in {\bf Z}}
((m+b)!)^{-1}\sum_{l=0}^m{{m+b}\choose m}
{m\choose l}$$
$$\{ (b+m-l)({\hat P}_{u^{b+m-l-1},w(u,\omega )^l}[(\partial ^{(m+b)}f/
\partial u^b\partial x^m)
(u,\xi (u,\omega ))\circ (I^{\otimes b}\otimes 
a^{\otimes (m-l)}\otimes E^{\otimes l})])|_{u=t}+$$
$$l({\hat P}_{u^{b+m-l},w(u,\omega )^{l-1}}[(\partial ^{(m+b)}f/
\partial u^b\partial x^m)
(u,\xi (u,\omega ))\circ (I^{\otimes b}\otimes 
a^{\otimes (m-l)}\otimes E^{\otimes (l-1)})]E{w'}_u(u,\omega ))
|_{u=t} \} .$$}
\par {\bf Proof.} In view of Theorem I.4.8 \cite{lust1} and
Proposition 3.11 there exists a generating operator of an evolution family.
From Lemma I.2.3 and Formula $I.4.8.(ii)$ \cite{lust1}
it follows the statement of this Corollary.
\par {\bf Remark.} If $f(t,x)$ satisfies conditions either of \S I.4.6
or of \S I.4.7, then Formula $3.12.(1)$ takes 
simpler forms, since the corresponding terms vanish. \\
\par The author is sincerely grateful to I.V. Volovich for
his interest to this work and fruitful discussions.

\end{document}